\documentclass[10pt,a4paper]{article}
\usepackage[utf8]{inputenc}
\usepackage[english]{babel}
\usepackage{amsmath}
\usepackage{amsfonts}
\usepackage{amssymb}
\usepackage{amsthm}
\usepackage{tikz-cd}
\usepackage{mathrsfs}
\usepackage{mathtools}
\usepackage{todonotes}

\newcommand{\V}{\mathbb{V}}
\newcommand{\Z}{\mathbb{Z}}
\newcommand{\N}{\mathbb{N}}

\newcommand{\im}{\mathop{\textup{im}}}

\newcommand{\trcl}{\mathop{\textup{tr cl}}}

\newcommand{\crit}{\mathop{\textup{crit}}}
\newcommand{\E}{\mathscr{E}}
\newcommand{\F}{\mathfrak{F}}

\renewcommand{\r}{\rightarrow}
\newcommand{\stab}{\textup{Stab}}
\newcommand{\rel}{\mathfrak{r}}
\renewcommand{\S}{\mathfrak{S}}

\author{Dianthe Basak \thanks{This project has received funding from the European Union’s Horizon 2020 research
		programme under the Marie Sklodowska-Curie grant agreement 945322. \includegraphics[scale=0.02]{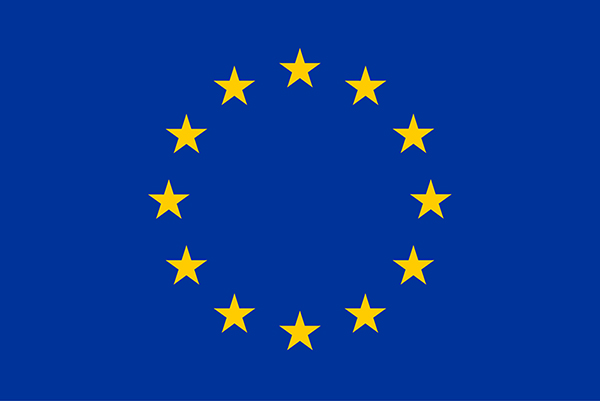}}}
\title{Monoidal Symmetric Models}

\newtheorem*{theorem*}{Theorem}
\newtheorem{theorem}{Theorem}[section]

\newtheorem{cor}{Corollary}[theorem]
\newtheorem{claim}[theorem]{Claim}
\newtheorem{lemma}[theorem]{Lemma}
\newtheorem{prop}[theorem]{Proposition}
\newtheorem{defn}[theorem]{Definition}

\begin{document}
\maketitle

\begin{abstract}
We develop a new method of interpreting large cardinal axioms as giving rise to topological symmetries of the universe of sets, similar to the construction of Fraenkel-Mostowski-Specker models. This allows us to define a ``symmetric" inner model construction. In this vein, we use Fraenkel-Mostowski-Specker type stabiliser arguments to deduce the Kunen inconsistency, as well as resolving an open problem of Rogers \cite{Rogers2023toposesof}.
\end{abstract}

\section{Introduction}

Since at least the work of Blass, Scedrov \cite{blassscedrov}, Freyd \cite{FREYD198749}, and contemporaries, it has been known that a topos theoretic viewpoint exists on most known methods for building models of ZFC. For instance, Blass and Scedrov interpret forcing as internal truth in the Grothendieck topos of canonical sheaves on a Boolean algebra, while Freyd's representation theorem exhibits topoi as exponential subvarieties inside a forcing topos constructed over the topos of $G$-sets for a topological group $G$ (on the set theoretic side, these are the names of a symmetric extension). As noted by Blass and Scedrov, since the inclusion functors of exponential Boolean subvarieties are logical functors, the entire logical content of a full Grothendieck topos is, roughly, always exhibited by a symmetric extension.

This is, however, not the end of the story. The correspondence above is by no means a complete translation of set theoretic methods into topos theoretic ones. Being a Grothendieck topos is itself a second-order notion, and thus inner models are completely missing from the symmetric extensions picture given above. Nevertheless, there are well-known topos theoretic methods to recover ``all" models of set theory, in the sense that they enjoy completeness theorems over ZF (though they are usually stated in intuitionistic set theory, these theorems essentially follow from Awodey Butz et al. \cite{AWODEY2014428}).

In the vein of constructing new Grothendieck topoi, one can ask which sites \cite{MacLane1994} give rise to new topoi that are of interest to set theorists. One such choice is that given by a topological monoid, giving rise to the topos of continuous actions of said monoid. Inspired by this, this article hopes to provide a few new ways of constructing inner models by using continuous actions of monoids on the universe of sets.

\section{Preliminaries}

\subsection{Fraenkel-Mostowski-Specker Models}

For simplicity, we work in ZFA, where the axiom of extensionality is replaced with
\[
[\exists z(z \in x) \wedge \forall z(z \in x \leftrightarrow z \in y)] \r x = y.
\]
and we add the axiom $\exists A \forall x (\forall y(y \notin x) \r x \in a)$. This permits the universe to be decomposed into the set of atoms $A$ (which have no elements and are not beholden to extensionality) and the rest. Thus the universe is built up by the strata $V_0(A) = A$, $V_{\alpha+1}(A) = A \cup P(V_\alpha(A))$, and $V_\lambda(A) = \bigcup_{\alpha < \lambda} V_\alpha(A)$, and $\V(A) = \bigcup_\alpha V_\alpha(A)$.

Any action of a group $G$ on $A$ extends recursively to the full universe $\V(A)$:
\[
g(S) = \{gx : x \in S\}
\]
If $(G, \tau)$ is a topological group whose topology $\tau$ is generated by a system of open subgroups of the identity, then we can ask about the continuity of the action of $G$ on $\V(A)$, i.e. the continuity of
\[
G \times \V(A) \to \V(A)
\]
where $\V(A)$ has the discrete topology. This simplifies to the well-known stabiliser condition: for each $x \in \V(A)$ we must have $\stab(x)$ open. But this is clearly not always true. For instance, given an action of $\Z$ with the topology of arithmetic progressions with positive common difference, if $x_n$ is a sequence of atoms with stabiliser $n\Z$, then the sequence $(x_i)_i$ has stabiliser $\bigcap_{n > 0} n\Z = \{0\}$ which is not open.

Hence we pass to the largest transitive submodel on which the action is continuous, namely
\[
\S(G, \tau; A) = \{x \in \V(A) : \forall y \in \trcl\{x\} (\stab(y) \text{ is open})\}
\]
The insight that long sequences $(x_i)_i$ have non-open stabilisers (are ``unstable") now manifests as the absence of sufficiently many sequences from $\V(A)$, hence failures of the axiom of choice. For instance, if $A$ is a disjoint union of the orbits of $x_n$, each with stabiliser $n\Z$ as before, then by a result of Levy \cite{Levy1962}
\[
\S(\Z ; A) \vDash ZF + |A|\text{ is incomparable with }\aleph_0
\]

\subsection{Monoid Actions}

The utility of atoms in the previous approach was to allow actions of groups, since isomorphisms from a transitive set with atoms to itself are completely determined by their action on the atoms. In particular, there is no non-identity isomorphism of the traditional universe of sets without atoms, and thus no action of a non-identity group on $\V = \V(\{\})$.

There are, however, many monoids that act on $\V$ with varying degrees of structural preservation. For instance, Hamkins et al. \cite{hamkins} establish that the map defined recursively by
\[
j(x) = \{j(y) : y \in x\} \cup \{\{0, x\} \}
\]
is an embedding of structures from $(\V, \in) \to (\V, \in)$ that is not the identity. Furthermore, under the assumption of the existence of a measurable cardinal, we have the existence of embeddings $j : \V \to \V$ which preserve and reflect all $\Delta_1$ predicates. In general we will write $j : \V \to_n \V$ to denote that $j$ preserves and reflects $\Sigma_n$ (equivalently preserves and reflects $\Delta_{n+1}$) predicates. It is a classical theorem that preservation and reflection of $\Delta_2$ predicates implies all higher $\Delta_n$ preservation and reflection. We would like to reproduce the construction of Fraenkel Mostowski models above with monoids of elementary embeddings.

To this end, we assume that $M$ is a monoid.
\begin{defn}
	A \textbf{principal $M$-set} is a quotient $M/R$ where $R$ is a \textbf{left congruence}, i.e. an equivalence relation on $M$ such that $xRy \r ax R ay$. Any $M$-set generated by a single element is of this form, hence the word principal.
\end{defn}

\begin{defn}
	Suppose $M$ is a topological monoid (i.e. the multiplication is continuous). A \textbf{continuous $M$-set} is a discrete set $X$ with an action of $M$ whose action map is continuous, i.e. such that for $x, y \in X$, there is an open set of $m \in M$ such that $mx = y$. This is equivalent to saying that the \textbf{stabiliser relation $\rel_x$}, defined as $\{(m, m') : mx = m'x\}$ has open classes. A principal $M$-set $M/R$ is continuous iff $Rm^{-1} = \{(n, n') : nm R n'm\}$ has open classes for all $m \in M$, iff $R$ has open classes iff $R$ is open in $M \times M$.
\end{defn}

\begin{defn}
	\textbf{A map of $M$-sets} $f : X \to Y$ is a function which commutes with the action of $M$.
\end{defn}

\begin{defn}
	\textbf{A left powder topology on $M$} is a topology which is generated by the continuous principal left $M$-sets. That is, a neighbourhood basis for every $x$ is given by sets of the form $\{y : yRx\}$ where $R$ is an open left congruence.
	
	Furthermore, suppose we are given a set $P$ of left congruences such that
	\begin{itemize}
		\item $R_1, R_2 \in P \r R_1 \cap R_2 \in P$
		\item $R_1 \subseteq R_2$, $R_1 \in P \r R_2 \in P$
		\item $R \in P \r Rm^{-1} \in P$
		\item If for every $x$ in $M$ there is $S \in P$ so that $[x]_S \subseteq [x]_R$, then $R \in P$.
	\end{itemize}
	Then we may define the \textbf{topology on $M$ generated by $P$} by taking $\{y : yRx\}$ as a basis. This declares all relations in $P$ to be open, and these are all the open left congruences with respect to this topology by the fourth condition above.
\end{defn}

In analogy with the group case,
\begin{defn}
	Suppose $(M, \tau)$ is a monoid with a left powder topology and an action on $\V(A)$, and define \textbf{the symmetric core}
	\[
	\S(M, \tau ; A) = \{x \in \V(A) : \forall y \in \trcl\{x\} (\rel_x \text{ is open})\}
	\]
	If $A = 0$ we will usually denote this $\S(M, \tau)$. If $\tau$ is the minimal topology in which the relations $\{\rel_x : x \in P\}$ are open, we will by abuse of notation write $\S(M, P)$ to mean $\S(M, \tau)$.
\end{defn}

\section{Examples}

Similar to the proof in \cite{Kunen1973} one can use stabiliser arguments to prove the failure of choice in many $\S(M, P)$ models.

\begin{theorem}
	Let $P := Ord^\omega$, and let $M$ be the monoid of all ultrapower embeddings by countably complete ultrafilters on ordinals in $\V$. Suppose there are $\omega_1$ many measurable cardinals with supremum $\eta$. Then
	\[
	\S(M, P) \models\text{there is no well-ordering of }[\eta]^\omega
	\]
\end{theorem}
\begin{proof}
Let $F : \alpha \to [\eta]^\omega$ be some purported well-ordering. Let $p = \{\kappa_n\}_{n < \omega}$ be a parameter so that
\[
[id]_{\rel_p} \subseteq [id]_{\rel_F}.
\]
We claim that
\begin{itemize}
	\item There is $s = \{\lambda_n\}_n \in [\eta]^\omega$ so that for every ordinal $\beta$, $[id]_{\rel_p \cap \rel_\beta} \not\subseteq [id]_{\rel_s}$.
	\item There is $\beta < \alpha$ so that $\rel_F \cap \rel_\beta \subseteq \rel_s$
\end{itemize}
These two facts are clearly in contradiction. The first follows by noting that for every ordinal $\gamma \in p$, there are only finitely many ordinals which support ultrafilters which move $\gamma$ (this is Kunen's classical lemma \cite{Kunen1973}), hence after countably many measurables, all remaining measurables $< \eta$ must support only ultrafilters that fix all of $q$. Let $s$ be any such sufficiently large countable sequence.

Then for any ordinal $\beta$, after finitely many elements of $s$, all remaining elements have no ultrafilters that move $\beta$. Hence such embeddings are $= id$ on $q \cup \{\beta\}$ but differ from the identity on $s$, as needed.

The second fact follows from noting that $s$ is in the image of $F$ by some ordinal $\beta$.
\end{proof}
Note that the above $\S(M, P)$ is an admissible structure containing $Ord^\omega$, but not necessarily a model of full ZF. However, there are situations where we do get ZF:

\begin{theorem}
Let $M$ be any OD monoid and let $P \supseteq \{V_\alpha\}_{\alpha \in Ord}$ be OD. Then $\S(M, P)$ is an inner model of ZF.
\end{theorem}
\begin{proof}
If $b$ is $\Delta_1$-definable from $a$, $\rel_a \subseteq \rel_b$. Hence if $\rel_a$ is open so is $\rel_b$. Thus $\S(M, P)$ is closed under Godel operations. Almost universality follows from the ordinal definability of $V_\alpha \cap \S$, since if this set is OD inside $V_\mu$,
\[
\rel_{V_\mu} \subseteq \rel_{V_\alpha \cap \S}
\]
making the latter open.
\end{proof}

\section{Non-discrete Actions and Kunen's Theorem}

The Kunen inconsistency is the celebrated result of Kunen \cite{Kunen_1971} that
\begin{theorem*}[Morse-Kelly set theory]
There is no fully elementary embedding $\V \to \V$.
\end{theorem*}

The core step of Woodin's proof of the Kunen inconsistency is the observation that a stationary subset of $cof_\omega(\lambda^+)$ intersects the image of any elementary embedding that is continuous at $\lambda^+$. This is because the image of an $\omega$-club in $\lambda^+$ under an elementary $j$ is still an $\omega$-club of $\lambda^+$. Similarly:

\begin{lemma}
If $j, k : \V \to_0 \V$ and $cof(\delta) > \omega$, $j(\delta) = k(\delta) = \delta$, then $\{\alpha < \delta : j(\alpha) = k(\alpha)\}$ is an $\omega$-club.
\end{lemma}
\begin{proof}
$\omega$-closure is clear. Suppose $\alpha = \alpha_0 < \delta$. Then we define $\alpha_{n+1} = \max\{j(\alpha_n), k(\alpha_n)\}$. Then the limit of the $\alpha_n$ is a fixed point of both $j$ and $k$, and is $\geq \alpha$.
\end{proof}

This leads to a version of Woodin's observation that is not ``lopsided":

\begin{lemma}
Suppose $f : cof_\omega(\lambda^+) \to X$ has $\omega$-stationary fibres, and $j, k : \V \to_0 \V$ satisfy $j(\lambda^+) = k(\lambda^+) = \lambda^+$. If $j(f) = k(f)$, then $j|_X = k|_X$.
\end{lemma}
\begin{proof}
Let $x \in X$, and choose $\alpha \in f^{-1}(x)$ at which $j(\alpha) = k(\alpha)$ (since $f^{-1}(\alpha)$ is $\omega$-stationary). Then
\[
j(x) = j(f(\alpha)) = j(f)(j(\alpha)) = k(f)(k(\alpha)) = k(x).
\]
\end{proof}

%

Using this we can rephrase the original proof of the inconsistency in topological terms:

\begin{prop}[Choice]
Let $\E$ be the monoid of all $\Delta_0$-embeddings $H(\lambda^+) \to H(\lambda^+)$, and let this monoid act on $P_{\lambda^{++}}(H(\lambda^+))$ by
\[
j^+(A) = \bigcup_{x \in H(\lambda^+)} j(A \cap x)
\]
Then the only topology on $\E$ with respect to which this action is continuous is the discrete topology.
\end{prop}
\begin{proof}
Suppose the embedding $j \in \E$ is not an isolated point. Then consider the open congruence $\rel_f$ where $f$ is some code for a stationary partition of $cof_\omega(\lambda^+)$ into $|V_\lambda|$ many stationary sets (which can be made a $\lambda^+$-size subset of $H(\lambda^+))$). Then any $k$ such that $j \rel_f k$ is equal to $j$ on $V_\lambda$. We compute first, for any subset $A \subseteq \lambda$, that
\[
\bigcup_{\alpha_n < \lambda} j(A \cap \alpha_n) = \bigcup_{\alpha_n < \lambda} k(A \cap \alpha_n)
\]
for any countable cofinal sequence $\alpha_n$, whence $j(A) = k(A)$. Thus, via coding, $j|_{H(\lambda^+)} = k|_{H(\lambda^+)}$. However then $j = k$, hence $\rel_f$ isolates $j$, a contradiction.
\end{proof}

\begin{prop}
If there is a $\Sigma^1_1$-elementary $j, j^* : H(\lambda^+) \to H(\lambda^+)$ that is $\Sigma^1_1$ (i.e. the class part, $j^*$, is $\Sigma_1$), then the action of $\E$ on $P_{\lambda^{++}}(H(\lambda^+))$, topologised by the least topology making the action continuous, is not discrete.
\end{prop}
\begin{proof}
If $p \in P_{\lambda^{++}}(H(\lambda^+))$ is any parameter, consider the formula
\[
\exists k : H(\lambda^+) \xrightarrow{\Delta_0} H(\lambda^+), k(k|_{V_\lambda}) = j|_{V_\lambda}, \bigcup_{x \in H(\lambda^+)} k(p \cap x) = \bigcup_{x \in H(\lambda^+)} j(p \cap x)
\]
The translation of this by $j^*$ is witnessed by $j$ itself, since
\[
\bigcup_{x \in H(\lambda^+)} j(j^*(p) \cap x) = \bigcup_{x \in H(\lambda^+)} j^*(j|_{H(\lambda^+)})(j^*(p) \cap x)
\]
This is because, while $j$ is \textbf{not} cofinal (indeed, by results of \cite{GoldbergPeriodicity}, it cannot be), we can select an enumeration $f : \lambda^+ \to p$ to get
\[
\bigcup_{\alpha < \lambda^+} j(\im j^*(f|_\alpha)) = \bigcup_{\alpha < \lambda^+} j^*(j|_{H(\lambda^+)})(\im j^*(f|_\alpha))
\]
whence we get that both sides are equal to
\[
j^*(j^*(\im f)) = j^*(j^*(p))
\]
essentially since $j$ is cofinal on ordinals upto $\lambda^+$. Since $j^*$ is $\Sigma_1$, we get a $k$ witnessing the above formula.

Thus the natural action of $k$ on $P_{\lambda^{++}}(H(\lambda^+))$ gives an element of the open neighbourhood $[j]_{\rel_p}$ which is not $j$ (since $\crit(k) < \crit(j)$). Thus no $\rel_p$, i.e. no open neighbourhood, isolates $j$.
\end{proof}
\begin{cor}[Choice]
There is no $\Sigma^1_1$ elementary $j, j^* : H(\lambda^+) \to H(\lambda^+)$.
\end{cor}
In the above proof $j^*$ need not be (indeed, never is) $j^+$. Rather, one can apply the above to situations where one has excess information in $j^*$ that one can restrict to an embedding of $H(\lambda^+)$.
\begin{cor}[Choice]
There is no elementary embedding $V_{\lambda+2} \to V_{\lambda+2}$.
\end{cor}
\begin{proof}
Let $X \subseteq H(\lambda^+)$, then one can code $X$ by
\[
\{c \in P(\lambda) : \text{ the Mostowski collapse of }c \text{ lies in }X\} \in V_{\lambda+2}
\]
The elementarity of $j$ then gives a $\Sigma^1_1$ elementary embedding $H(\lambda^+) \to H(\lambda^+)$.
\end{proof}

Notice that the core observation about stationarity did not need any strong embedding assumptions. Indeed the place we used the elementarity of $j$ was really that $j$ had an ``infinitesimal perturbation" $j(j)$ which was at once ``close" ($j(j)(j(f)) = j(j(f))$) to and ``far" (different critical points) from it. We note, with a view to potential applications, the following basic fact.

\begin{prop}
Let $M$, $N$ be models of ZF$^-$, $j, k : M \to_0 N$, and $\F$ be a filter on $X \in M$ ($\F$ is not necessarily in $M$). Suppose $j$ and $k$ are equal on an $\F$-large set ($j_*\F = k_*\F$) and $f : X \to Y$ is a function with $\F$-positive fibres. Then $j(f) = k(f)$ implies $j|_Y = k|_Y$ (in the above notation, $\rel_{\F} \cap \rel_f \subseteq \bigcap_{y \in Y} \rel_y$).
\end{prop}

%

\section{A Problem of Rogers}

In his paper \cite{Rogers2023toposesof}, Rogers poses the question of whether there exists a left powder monoid that is not a right powder monoid. We recall the definition of these terms in the equivalent form given by Rogers:

\begin{defn}
	A \textbf{left powder monoid} is a topological monoid $(M, \tau)$ which is T0 and such that $\tau$ has a basis of clopen sets $U$ for which
	\[
	I^U_p = \{q : \{r : rq \in U\} = \{r : rp \in U\}\} \in \tau
	\]
	for every $p$. Note that this is just $[U]_{\rel_p}$ for the appropriate action, but we reserve the latter notation only for actions on the universe of sets. A \textbf{right powder monoid} is a topological monoid $(M, \tau)$ such that $(M^{op}, \tau)$ is a left powder monoid. For ease of notation, we will call a left powder monoid that is not right powder (with the same topology) a \textbf{chiral monoid}.
\end{defn}

This definition of being left powder expresses in essence that the topology on $M$ is derived from its action on some discrete set (and that this set can be taken to be $\tau$, see \cite{Rogers2023toposesof}). Note that the left action on $\tau$ is induced by $M$'s right action on itself.
%
%

\begin{defn}
	Let $M = \{j : V_\lambda \to_1 V_\lambda\}$ where $\lambda$ is limit. \textbf{The canonical topology} $\tau$ on $M$ is the one that makes its action on $V_\lambda$ continuous, i.e. the basic opens around $j$ are
	\[
	[j]_{\rel_x} = \{k \in M : k(x) = j(x)\}
	\]
	for $x \in V_\lambda$. We will sometimes abuse notation and denote these $[j]_x$.
\end{defn}

\begin{lemma}
	The topology $(M, \tau)$ is $T_0$. Since it is generated by clopens, it is totally disconnected Hausdorff.
\end{lemma}
\begin{proof}
	If $j, k \in M$ are inseparable by any open sets, then for every $x$, $j(x) = k(x)$ (since $j \in [k]_x$). But then they are equal as functions.
\end{proof}

\subsection{Chirality}

\begin{claim}
	$(M, \tau)$ is a left powder monoid.
\end{claim}
\begin{proof}
	As noted before, this follows essentially from the definition of $M$, as it is topologised by its action on $V_\lambda$. For completeness, we include a proof.
	
	Let $U = [j]_x$,
	\[
	I^U_p = \{q : \{r : rq(x) = j(x)\} = \{r : rp(x) = j(x)\}\}.
	\]
	Suppose $q_0 \in I^U_p$, and consider $q \in [q_0]_x$, then for any $r$,
	\[
	rq(x) = j(x) \leftrightarrow rq_0(x) = j(x) \leftrightarrow rp(x) = j(x)
	\]
	and hence $q \in I^U_p$, i.e. $[q_0]_x \subseteq I^U_p$. Since $q_0$ was arbitrary, $I^U_p$ is open.
\end{proof}

The crucial tool to show that $M$ is chiral, and indeed the only part of the proof that exploits the special structure of elementary embeddings, is the following lemma:

\begin{lemma}[Laver]
	\label{laver}
	Suppose $j : V_{\lambda+1} \to V_{\lambda+1}$ is elementary. We consider its restriction to $V_\lambda$ which lies in $M$. For any parameter $s, t \subseteq V_\lambda$, and any ordinal $\alpha < \crit j$, there is $r \in M$ so that
	\[
	r^+(r) = j|_{V_\lambda}, r^+(s) = j(s), t \in \im r^+, \crit(r) > \alpha.
	\]
\end{lemma}
\begin{proof}
	Fix $s, t \subseteq V_\lambda$, $\alpha < \crit j$. We apply $j$ to the above statement, i.e.
	\[
	\exists r \in M (r^+(r) = j(j|_{V_\lambda}), r^+(j(s)) = j(j(s)), j(t) \in \im r^+, \crit(r) > j(\alpha) = \alpha.
	\]
	This formula is witnessed by $r = j|_{V_\lambda}$ (note that $M$ has not changed in the formula since ``$r : V_\lambda \to V_\lambda$ is elementary" is definable in $V_{\lambda+1}$). Hence the original formula must also have a witness by the elementarity of $j$.
\end{proof}
By a very similar analysis, we obtain:
\begin{claim}
	$(M, \tau)$ is not a right powder monoid.
\end{claim}
\begin{proof}
	Suppose $j : V_{\lambda+1} \to V_{\lambda+1}$ with critical point $\kappa$. Assume for contradiction that $\tau$ has a basis for which
	\[
	I^U_p = \{q : \{r : qr \in U\} = \{r : pr \in U\}\}
	\]
	is open for each $U$ in the basis, $p \in M$.
	
	Set $J = j \circ j(j|_{V_\lambda})$ with critical point $\kappa$. We claim that for \textbf{any} open $J \in U \subseteq [J]_\kappa$,
	\[
	I^U_{j(j)} = \{q : \{r : qr \in U\} = \{r : j(j)r \in U\}\}
	\]
	contains no open neighbourhood of $j(j)$ despite containing $j(j)$ by definition. Thus $I^U_{j(j)}$ is not open, hence the basis of $\tau$ above does not cover $J$, a contradiction.
	
	To prove the claim we prove two lemmas.
	\begin{lemma}
	\[
	\forall r (j(j)r(\kappa) \neq J(\kappa))
	\]
	\end{lemma}
	\begin{proof}
	If $j(j)(r(\kappa)) = j(j(j)(\kappa))$, then $j(j(j)(\kappa)) \in \im j(j)$, whence $j(j)(\kappa) \in \im j$. But $j(j)(\kappa) = \kappa$ is not in the image of $j$, a contradiction.
	\end{proof}
	\noindent Then $\{r : j(j) r \in U\} \subseteq \{r : j(j)r \in [J]_\kappa\} = \{\}$, whence
	\[
	I^U_{j(j)} = \{q : \forall r (qr \notin U)\}
	\]
	Now let $J \in [J]_v \subseteq U$ be any neighbourhood smaller than $U$, where $v$ is some element of $V_\lambda$ (these are the basic opens).
	\begin{lemma}
	\[
	\forall x \exists q : V_\lambda \to V_\lambda (q(x) = j(j)(x) \wedge \exists r : V_\lambda \to V_\lambda (qr(v) = J(v)))
	\]
	In other words, there is an element $q \in [j(j)]_x$ so that $qr \in [J]_v \subseteq U$ for some $r$, hence $q \notin I^U_{j(j)}$.
	\end{lemma}
	\begin{proof}
	Fix $x$. Push the formula
	\[
	\exists q (q(x) = j(j)(x) \wedge \exists r (qr(v) = J(v)))
	\]
	forward by $j(j)$ to get
	\[
	\exists q (q(j(j)(x)) = j(j)(j(j)(x)) \wedge \exists r (qr(j(j)(v)) = j(j)(J(v))))
	\]
	It is not hard to see that $q = j(j)$, $r = j$ witness this formula, hence by the elementarity of $j(j)$ we are done.
	\end{proof}
	\noindent These two lemmas establish that $I^U_{j(j)}$ is not open.
\end{proof}

Upon inspecting the above proof, the reader may note a general framework for obtaining a chiral monoid. In general, the above result can be achieved by arranging a monoid $M$ with topology derived from its action on $X$, with elements $A, B \in M$ so that
\begin{itemize}
	\item $\neg\exists r (B\circ r = A)$
	\item $\forall x, v \in X^{<\omega} \exists q, r \in M (B(x) = q(x), q \circ r(v) = A(v))$.
\end{itemize}
Equivalently, the map $(q, r) \mapsto (q, q \circ r)$ does not have closed image. With the knowledge that such conditions can be achieved from the assumption of embeddings $V_{\lambda+1} \to V_{\lambda+1}$, it is not hard to see that simpler monoids also may satisfy them (for instance, $\N^\N$).

However, in light of lemma \ref{laver} the endomorphism monoids of ranks were a natural first guess. This forms part of a long history of work by Dehornoy \cite{Dehornoy2010}, Laver \cite{LAVER1995334}, Dougherty \cite{DOUGHERTY1993211}, and others deriving combinatorial results using strong set theoretic assumptions that could later be eliminated from the proof.

\subsection{Completeness}

Let $\mathscr{R}$ be the set of open left congruences on $M$, and consider the map
\[
M \to L := \varprojlim_{R \in \mathscr{R}} M/R
\]
where the limit is over the subcategory of continuous $M$-sets on the objects $M/r$ with quotient maps $M/R \to M/R'$ whenever $R \subseteq R'$. The map in question takes $m \in M$ to the constant sequence $([m]_R)_{R \in \mathscr{R}}$.

\begin{defn}
	A topological monoid $(M, \tau)$ is \textbf{left complete} if this map is an isomorphism of topological monoids.
\end{defn}

\begin{claim}
	$(M, \tau)$ is left complete.
\end{claim}
\begin{proof}
	Note that since we obtained the topology on $M$ from an action, this map is automatically injective: if $j \neq k$, then choosing any $x \in V_\lambda$ on which they differ, we note that $\rel_x$ is an open left congruence, and $[j]_{\rel_x} \neq [k]_{\rel_x}$, hence the images of $j$ and $k$ in $L$ are different.
	
	Let $[\alpha_R]_R$ be any sequence representing an element of $L$. Define a map $\alpha : V_\lambda \to V_\lambda$ by
	\[
	\alpha(x) = \alpha_{\rel_x}(x)
	\]
	which is well-defined. We claim this map is also $\Delta_0$ elementary and cofinal. This will show that $\alpha \in M$. To show cofinality, note
	\[
	\alpha(V_\alpha) \cap V_\alpha = V_\alpha
	\]
	since the codomains of all the $\alpha_r$ are $V_\lambda$. Note that since $\alpha_R$ is a sequence in $L$, if $R \subseteq R'$, $[\alpha_R]_{R'} = [\alpha_{R'}]_{R'}$, thus in particular
	\[
	\alpha_{\rel_x \cap \rel_y}(x) = \alpha_{\rel_x}(x).
	\]
	Suppose $\phi(x, \bar{y})$ is a $\Delta_0$ formula with free variables shown, then
	\[
	V_\lambda \models \phi(x, \bar{y}) \leftrightarrow V_\lambda \models \phi(\alpha_{\rel_x \cap \rel_{\bar{y}}}(x), \alpha_{\rel_x \cap \rel_{\bar{y}}}(\bar{y}))
	\]
	which is to say $V_\lambda \models \phi(\alpha_{\rel_x}(x), \alpha_{\rel_{\bar{y}}}(\bar{y}))$, whence $\alpha$ is $\Delta_0$ elementary. Hence $\alpha \in M$ is a preimage of $[\alpha_R]_R$ under the canonical map. Thus $M \simeq L$ as monoids. Since $M$ is powder it is an isomorphism of topological monoids.
\end{proof}

This resolves a second question posed by Rogers:
\begin{cor}
There is a left complete monoid which is not right powder (hence not right complete).
\end{cor}

\nocite{*}
\bibliographystyle{alpha}
\bibliography{article1_refs}
\end{document}